\documentclass[12pt]{article} 
\usepackage{amssymb}
\newtheorem{thm}{Theorem}
\newtheorem{lem}{Lemma} 

\begin{document}

\title{The minimum forcing number of perfect matchings in the hypercube}

\author{Ajit A. Diwan  \\ 
Department of Computer Science and Engineering,\\
Indian Institute of Technology Bombay, Mumbai 400076, India.}

\maketitle

\textbf{Abstract}
Let $M$ be a perfect matching in a graph. A subset $S$ of $M$ is said to be 
a forcing set of $M$, if $M$ is the only perfect matching in the graph that
contains $S$. The minimum size of a forcing set of $M$ is called the forcing
number of $M$.  Pachter and Kim [Discrete Math. 190 (1998) 287--294] conjectured 
that the forcing number of every perfect matching in the $n$-dimensional hypercube
is at least $2^{n-2}$, for all $n \ge 2$. Riddle [Discrete Math. 245 (2002) 283-292] proved 
this for even $n$. We show that the conjecture holds for all $n \ge 2$. The proof is 
based on simple linear algebra.

\section{Introduction}

Let $M$ be a perfect matching in a graph $G$. A subset $S$ of edges in $M$ 
is said to be a \emph{forcing set} of $M$, if $M$ is the only perfect matching
in $G$ that contains $S$.  The minimum size of a forcing set of $M$ is called
the forcing number of $M$.

Forcing sets and forcing numbers of perfect matchings were first studied by
Harary et al.~\cite{HKZ}, and have applications in chemistry~\cite{CC}. Pachter and
Kim~\cite{PK} found tight upper and lower bounds on the forcing number of
perfect matchings in a 2-dimensional square grid. They also conjectured that every
perfect matching in the $n$-dimensional hypercube has forcing number at least
$2^{n-2}$, for all $n \ge 2$. This was proved by Riddle~\cite{R} for even $n$. 
Adams et al.~\cite{AMM} showed that for all $n \ge 5$ and integers $i$ in the 
interval $[2^{n-2}, 2^{n-2}+2^{n-5}]$, there exists a perfect matching in the 
$n$-dimensional hypercube with forcing number $i$.

In this note, we complete the proof of Pachter and Kim's conjecture, and show that
it holds for all $n \ge 2$. Our proof is much simpler than the one in~\cite{R}, and uses
only basic linear algebra. We actually prove a slightly stronger statement. If $S$
is a forcing set of a perfect matching $M$ in a graph, then the graph obtained
by deleting all endvertices of edges in $S$ has a unique perfect matching $M\setminus S$.
Thus an upper bound on the order of an induced subgraph with a unique perfect
matching directly gives a lower bound on the forcing number of every perfect
matching in the graph. We show that any induced subgraph of the $n$-dimensional 
hypercube, with a unique perfect matching, has at most $2^{n-1}$ vertices, for $n \ge 2$. 
This implies the lower bound on the forcing number. The technique used is quite general, 
and can be applied to any bipartite graph.

\section{Hypercube}

We describe our technique in general, and then apply it to the specific case of
the hypercube. Let $G$ be a bipartite graph, with a bipartition $X,Y$ of the 
vertex set. Let $X = \{x_1,\ldots,x_n\}$ and $Y = \{y_1,\ldots,y_m\}$.  
The weighted bipartite adjacency
matrix of $G$ is an $n \times m$ matrix $W$, with $W_{ij} = 0$ if $x_i$
is not adjacent to $y_j$, and  $W_{ij} = w$ otherwise. Here, $w$ is an
indeterminate used to indicate the presence of an edge in the graph.

Let $H$ be an induced subgraph of $G$ with a unique perfect matching. 
The rows and columns of $W$ that correspond to vertices in $H$
form a square submatrix $W_H$ of $W$. Since $H$ has a unique perfect
matching, the determinant of $W_H$ must be $\pm w^{|H|/2}$, since
only one term in the expansion of the determinant is non-zero. This implies
that if we replace each occurrence of $w$ in $W$ by an arbitrary \emph{non-zero}
element from some field $\mathcal{F}$, the resulting matrix must have rank at 
least $|H|/2$ over $\mathcal{F}$.  Thus, if we choose a field $\mathcal{F}$, and 
appropriate non-zero values from the field as the weights of the edges, then twice the rank 
of the resulting bipartite adjacency matrix is an upper bound on the order of any 
induced subgraph with a unique perfect matching. We now apply this technique to 
the hypercube.

The vertex set of the $n$-dimensional hypercube $Q_n$ is $\{0,1\}^n$, the set of all
sequences of length $n$ in which each element is either $0$ or $1$. Two sequences
are adjacent if they differ in exactly one position. Let $E_n$ denote the subset of
these sequences with an even number of $1$ elements, and $O_n$ the subset
with an odd number of $1$ elements. The hypercube is a bipartite graph with $E_n, O_n$
the two parts of the bipartition.  Construct the weighted bipartite adjacency matrix $W_n$
of $Q_n$ with rows indexed by elements of $E_n$ in lexicographic order, and columns
indexed by elements of $O_n$ in lexicographic order.  Then $W_n$ is a symmetric matrix
with a simple recursive structure. $W_1$ is the $1 \times 1$ matrix with $w$ as the only entry
and
\begin{eqnarray*}
W_{n+1} & = & \left[
\begin{array}{c|c}
W_n & wI_{2^{n-1}} \\
\hline
wI_{2^{n-1}} & W_n
\end{array}
\right].
\end{eqnarray*}

Here, $I_k$ denotes the $k \times k$  identity matrix.

We use the recursive structure of $W_n$ to assign non-zero values from the field 
$Z_3$ of integers modulo $3$, to the $w$ elements in $W_n$, so that the rank of 
the resulting matrix is $2^{n-2}$ over $Z_3$. 

\begin{lem}
\label{lem1}
There exists a nonsingular matrix $A_n$, obtained by assigning non-zero values from 
the field $Z_3$ to the $w$ elements in $W_n$, such that $A_n^{-1}$ can 
also be obtained by assigning non-zero values to the $w$ elements in $W_n$.
\end{lem}

\noindent \textbf{Proof: }
We construct the matrix inductively. For $n = 1$, let $A_1$ be the $1 \times 1$
matrix containing the entry $1$.  Clearly, $A_1^{-1} = A_1$ satisfies the required
property. Let 
\begin{eqnarray*}
A_{n+1} & = &
\left[
\begin{array}{c|c}
2A_n & I_{2^{n-1}}\\
\hline
I_{2^{n-1}} & A_n^{-1}
\end{array}
\right]
.
\end{eqnarray*}

Then it is easy to verify that 
\begin{eqnarray*}
A_{n+1}^{-1} & = &
\left[
\begin{array}{c|c}
 A_n^{-1} & 2I_{2^{n-1}}\\
\hline
2I_{2^{n-1}} & 2A_n
\end{array}
\right].
\end{eqnarray*}

By induction,
and using the recursive structure of $W_{n+1}$, we conclude that both $A_{n+1}$
and $A_{n+1}^{-1}$ can be obtained from $W_{n+1}$ by assigning non-zero values
to the $w$ elements in $W_{n+1}$.
\hfill $\Box$

\begin{thm}
Any induced subgraph of $Q_n$ that has a unique perfect matching
contains at most $2^{n-1}$ vertices, for all $n \ge 2$.
\end{thm} 

\noindent \textbf{Proof: } Let $A_{n-1}$  be the matrix over $Z_3$ satisfying the 
properties in Lemma~\ref{lem1}. Let 
\begin{eqnarray*}
B_n & = & 
\left[
\begin{array}{c|c}
A_{n-1} & I_{2^{n-2}}\\
\hline
I_{2^{n-2}} & A_{n-1}^{-1}
\end{array}
\right].
\end{eqnarray*}

Now, it is easy to see that $B_n$ has rank $2^{n-2}$ over $Z_3$, and is obtained 
from $W_n$ by replacing all occurrences of $w$ in $W_n$ by non-zero values from
$Z_3$.  
\hfill $\Box$

Note that for even $n$, we get a simpler proof by considering the field $Z_2$.
In this case, if we replace all occurrences of $w$ in $W_n$ by $1$, we get a
matrix of rank $2^{n-2}$ over $Z_2$. However, for odd $n$, this gives a
nonsingular matrix which is its own inverse. Riddle's proof for even $n$ is
purely combinatorial. It would be interesting to see if the proof for odd $n$
can also be made purely combinatorial. It would also be interesting to see if
this technique can be used to find tight bounds on the forcing numbers of
perfect matchings in other bipartite graphs.

\end{document}